\documentclass[a4paper,11pt]{article}
\pdfoutput=1

\usepackage[utf8]{inputenc}
\usepackage{utf8}

\title{From Multisets to Sets in Homotopy Type Theory}
\author{Håkon Robbestad Gylterud}

\usepackage[hyphens]{url}
\usepackage{fancyvrb}

\usepackage{lmodern}
\usepackage{graphicx}

\usepackage{geometry}

\usepackage{longtable,booktabs}

\usepackage{unnumberedtotoc}

\usepackage[backend=bibtex,
            style=authoryear-icomp,
            hyperref=false,
            backref=false,
            maxcitenames=3
            ]{biblatex} 

\usepackage{amssymb,amsmath,amsthm}

\usepackage{stmaryrd}

\usepackage[all]{xy}
\input{diagxy}

\usepackage{chngcntr}
\counterwithout*{paragraph}{subsubsection}
\counterwithin*{paragraph}{section}

\theoremstyle{definition}

\newcounter{example}

\usepackage[explicit]{titlesec}

\renewcommand{\thesection}{\arabic{section}}

\renewcommand{\theparagraph}  {\thesection:\arabic{paragraph}}

\titleformat{\paragraph}[runin]{\normalfont\bfseries}{}{0em}{#1\ \theparagraph.}

\usepackage{xr}
\newcommand{\multisetsref}{\cite{multisets-arxiv}}

\usepackage{iterative-sets/notation}

\begin{document}

\maketitle

\begin{abstract}
We give a model of set theory based on multisets in homotopy type theory.
The equality of the model is the identity type. The underlying 
type of iterative sets can be formulated in Martin-Löf type theory, without
Higher Inductive Types (HITs), and is a sub-type of the underlying
type of Aczel's 1978 model of set theory in type theory. The Voevodsky
Univalence Axiom and mere set quotients (a mild kind of HITs) are used to
prove the axioms of constructive set theory for the model. We give an
equivalence to the model provided in Chapter 10 of “Homotopy Type Theory”
by the Univalent Foundations Program.
\end{abstract}

\newpage
\tableofcontents
\newpage

# Introduction

The first model of set theory in type theory is due to Aczel and models
the constructive set theory CZF\footcite{aczel1978}. The underlying type
of sets in this model is $W_U T$, the type of all well-founded trees
with branchings in a universe $U$ with decoding family
$T : U → {\Type}$. The interpretation of equality in this model allows
deduplication and permutation of subtrees — incorporating the intuition
that the order and multiplicity of elements of a set are irrelevant. If
we instead insist to interpret equality as the identity type and assume
the univalence axiom, the underlying type no longer models set theory,
but rather multiset theory.

In a related work^[\multisetsref] we explore 
the notion of iterative
multisets in type theory. Here we will specialise from the general
multisets to the set-like ones, where each element occurs at most once.
We will start by summarising the model of multisets and its relation to
Aczel’s model of CZF in type theory\footcite{aczel1978}.

Once the notion of a multiset is defined, it is natural to study the
hereditary subtype of multisets where each element occurs at most once.
These are in a certain sense the most natural representations of
iterative sets from a homotopy type theory point of view. These are
namely the multisets for which the element hood relation is hereditarily,
merely propositional (type level −1).

In this text we explore how this type models various axioms of
constructive set theory. We also show that it is equivalent to the
higher inductive type outlined in the book “Homotopy Type Theory”
\footcite{hottbook}.

## From multisets to sets – two ways

Assume for a moment that we would like to explain the notion of multiset
to someone who knows the notion of a set. Here are two similar attempts
at such an explanation.

### Variant A

> A multiset is a generalisation of sets in which an element may
> occur any number of times, not just at most once.

### Variant B

> A multiset is like a set with the extra information attached to
> each element about how many times it occurs in the multiset.

The two descriptions are almost identical, but critically different.
Variant A describes the multisets as a more general concept than
sets, in the sense that a set is a *special case* of a multiset,
namely those multisets in which each element occurs at most once.
Variant B, on the other hand, describes a multisets as sets with
some extra structure. In both cases the impression given is
that the multisets are, in an informal way, a *bigger concept* than
that of a set. Variant A says that sets are just some of the
multisets, while Variant B says that for each set there are numerous
ways to make a multiset from it.

In mathematics these two notions of “bigger”  correspond to the
notion of

 - sets being a subtype of multisets (Variant A)
 - sets being a quotient of multisets (Variant B)

If we were to turn the direction of explanations,
making them explanations of the concept of a set from the notion of
a multiset, the above suggests two distinct routes of constructing
a notion of sets from a notion of multisets. We either
identify the subtype of *set-like* multisets, or identify the
equivalence relation of multisets which *forgets the extra structure
of multiplicites*.

#### Example

The multiset $\{a,a,b\}$ would not be considered a set-like in the spirit
of Variant A, since $a$ occurs twice. On the other hand, Variant B would
contest that $\{a,a,b\}$ and $\{a,b,b\}$ would *represent the same*
set, since in both cases $a$ and $b$ are exactly the elements occurring
at least once.
\vspace{3mm}

If we denote by $M$ the multisets, and $V_A$ and $V_B$ denotes the
two possible notions of sets arising from it, following Variant A
and Variant B respectively, we can draw the following diagram of the
situation.

\begin{align*}
\xymatrix{
V_A \ar@{^{(}->}[r] & M \ar@{->>}[d] \\
  & V_B
}
\end{align*}

## Outline
In Section 2 and 3 we set up a bit of framework to work within. Starting
from Section 4 we will follow the path of Variant A.

Section 5 will take us through the basic lemmas about the type of
iterative sets we define in Section 4. These lemmas are applied in
Section 6 to give proofs that our model satisfies the axioms of Myhill’s
Constructive Set Theory.

In Section 7 we consider the problem of interpreting the two collection
axioms of Aczel’s CZF in our model.

In Section 8 we return to Variant B, which will make the relationship
with the approach taken in Chapter 10 of the book “Homotopy Type
Theory”\footcite{hottbook} clear.

## Notation

In what follows we will adhere to the following notation. Some of our
notation is similar to that of the book “Homotopy Type Theory”
\footcite{hottbook},
while some of it is inspired by the syntax of Agda.

 - Function application will be denoted by juxtaposition, as in $f␣a$.
   Also, application of functions equipped with extra structure, such as
   equivalences and embeddings, will be denoted by juxtaposition.
 - Quantifiiers, such as ∀,∃,∏ and ∑ bind weakly. For instance, \newline
   $∏_{x : ℳ} x∈a → ∑_{y : ℳ} {y∈B}∧{P␣x␣y}$ disambiguates to\newline
   $∏_{x : ℳ} \left(x∈a → ∑_{y : ℳ} {y∈B}∧{P␣x␣y}\right)$
 - The equality sign, $=$, denotes the identity type. We sometimes
   equip it with a subscript emphasising which type the elements
   belong to, as in $a =_A a'$.
 - Definitions are signified by $≔$.
 - Judgemental equalities are denoted by $≡$.
 - The notation $A : \Type$ denotes that $A$ is a type,
 - The notation $A : \Set$ denotes that $A$ is a type which is a mere set and
 - The notation $A : \Prop$ denotes that $A$ is a type which is a mere proposition.
 - The type $U$ is a universe, with decoding family $T : U → \Type$.
   We assume that this universe
     - is univalent,
     - contains the empty type, $0$,
     - is closed under Π-types,
     - is closed under Σ-types,
     - is closed under $+$-types,
     - is closed under $(-1)$-truncation and
     - is closed under taking quotients of mere sets
       by equivalence relations.
 - The notation $e_A$ denotes the empty function $e : 0 → A$. The subscript
   is dropped when inferable from the context.
 - The notation $ap␣f$ refers to the usual function $\ap␣f : a = a' → f␣a = f␣a'$
# Types and propositions

One of the main features of Martin-Löf’s type theory is the
interpretation of propositions as types\footcite{mltt}. The presence
of the identity
type gives the possibility of asking whether two proofs of a proposition
are equal. A type in which all elements are equal is called, in homotopy
type theory, a *mere proposition*.

The traditional interpretation of the existential quantifier in type
theory is by the sigma type $∑_{a : A} P␣a$. A proof of an existential
proposition is thus $a , p$ — a term $a : A$ of the quantified domain
paired with a proof $p : P␣a$ that the term has the correct property. It
is clear that since the existence is not necessarily unique, the type of
such proof need not be a mere proposition.

In homotopy type theory, one introduces a truncated existential
quantifier $∃(x:A) P␣a$, which is constructed from $∑_{a : A} P␣a$ by
adding identifications of all elements in to make it a mere proposition.
This gives the following introduction and elimination rule:^[In these
rules, written in an Agda-like notation, `Prop` refers to the type
of mere propositions.]

       a : A     b : B a
    ─────────────────────── ∃-intro
         [a,b] : ∃(a:A)(B a)

       x : ∃(a:A)(B a)   y : ∃(a:A)(B a)
    ────────────────────────────────── ∃-quot
         q : x = y

      P : ∃(a:A)(Ba) → Prop  
      p : (a : A) → (b : B a) → P [a,b]
    ──────────────────────────────────────────── ∃-elim
         ∃-elim P p : (x : ∃(a:A)(B a)) → P x

Clearly $∑_{a:A}(B␣a) → ∃(a:A)(Ba)$ holds for all $A$ and $B$. The
opposite implication holds if $∑_{a:A}(B␣a)$ is a mere proposition –
which is to say that the existence is unique, with a unique proof.

The situation is similar for disjunctions. Traditionally, disjunction is
interpreted as disjoint union in Martin-Löf’s type theory, while
homotopy type theory introduces a truncated variant. We will denote
disjoint unions by the operator $+$ and the truncated disjunction by the
operator $∨$.

# Models where equality is identity

In this section we define what an ∈-structure is, and give some basic
results on such structures in generality. We do this in order to adjust
our expectation for the concrete model which will be main focus of this
work.

## ∈-structures

#### Definition

An *∈-structure* is a pair $(ℳ,∈)$ where $ℳ : \Set$ is a mere set, and
$∈ : ℳ → ℳ → \Prop$.

#### Definition

For any ∈-structure $(ℳ,∈)$ and element $a : ℳ$, we define
$E␣a ≔ ∑_{x : ℳ} x∈a$.

#### Definition

An ∈-structure $(ℳ,∈)$ is called *$U$-like* if for each $a : ℳ$ the type
$E␣a$ is essentially $U$-small. That is, if each $a : ℳ$ the type
$E␣a$ has a code in $U$.

#### Remark

An ∈-structure is basically a mere set with a merely propositional, binary
predicate defined on it. The natural equality to consider for elements of
such a structure is the identity type. This is in contrast to the setoid
approach taken by Aczel.

“$U$-like” is meant to mimic the traditional terminology, “set-like”, used
in set theory.

## Translations of first-order logic into type theory

#### Definition

Given an ∈-structure $(ℳ,∈)$ define two translations of formulas of
first-order logic to type theory, $σ_{ℳ,∈}$ and $τ_{ℳ,∈}$, by recursion
on formulas. We let $\Ctx$ denote contexts (finite lists of
variables, $\Vars␣Γ$ denoting the variables in $Γ$) and let
 $\Fml␣Γ$ denote formulas in a given context $Γ : \Ctx$.

Starting with $σ_{ℳ,∈}$, leaving out the subscripts for ease of
reading:

\begin{align*}
  σ &: ∏_{Γ:\Ctx} \Fml␣Γ → (\Vars␣Γ → ℳ) → \Type \\
  σ &␣Γ␣(∀x\ φ)␣γ := ∏_{a:ℳ} σ␣(Γ.x) φ␣(γ.a) \\
  σ &␣Γ␣(∃x\ φ)␣γ := ∑_{a:ℳ} σ␣(Γ.x) φ␣(γ.a) \\
  σ &␣Γ␣(φ∧ψ)␣γ := σ␣Γ␣φ␣γ × σ␣Γ␣ψ␣γ \\
  σ &␣Γ␣(φ∨ψ)␣γ := σ␣Γ␣φ␣γ + σ␣Γ␣ψ␣γ \\
  σ &␣Γ␣(φ→ψ)␣γ := σ␣Γ␣φ␣γ → σ␣Γ␣ψ␣γ \\
  σ &␣Γ␣⊥␣γ := 0 \\
  σ &␣Γ␣⊤␣γ := 1 \\
  σ &␣Γ␣(x ∈ y)␣γ := γ␣x ∈ γ␣y\\
  σ &␣Γ␣(x = y)␣γ := γ␣x =_{ℳ} γ␣y
\end{align*}

\noindent where $Γ.x$ denotes the context $Γ$ extended with the
variable $x$, and $γ.a$ denotes the function $\Vars␣(Γ.x) → ℳ$
which maps $x$ to $a$.

We define $τ_{ℳ,∈}$ analogously, only difference being in
the clauses for ∨ and ∃:

\begin{align*}
  τ &: ∏_{Γ:\Ctx} \Fml␣Γ → (\Vars␣Γ → ℳ) → \Type \\
  τ &␣Γ␣(∀x\ φ)␣γ := ∏_{a:ℳ} τ␣(Γ.x) φ␣(γ.a) \\
  τ &␣Γ␣(∃x\ φ)␣γ := ∃(a:ℳ)\ τ␣(Γ.x) φ␣(γ.a) \\
  τ &␣Γ␣(φ∧ψ)␣γ := τ␣Γ␣φ␣γ × τ␣Γ␣ψ␣γ \\
  τ &␣Γ␣(φ∨ψ)␣γ := τ␣Γ␣φ␣γ ∨ τ␣Γ␣ψ␣γ \\
  τ &␣Γ␣(φ→ψ)␣γ := τ␣Γ␣φ␣γ → τ␣Γ␣ψ␣γ \\
  τ &␣Γ␣⊥␣γ := 0 \\
  τ &␣Γ␣⊤␣γ := 1 \\
  τ &␣Γ␣(x ∈ y)␣γ := γ␣x ∈ γ␣y\\
  τ &␣Γ␣(x = y)␣γ := γ␣x =_{ℳ} γ␣y
\end{align*}

#### Example

The axiom of union is translated differently by the two translations:

\begin{align*}
 \UNION ≔& ∀x␣∃u␣∀z\ z∈u ↔ ∃y∈x␣z∈y \\
 σ␣()␣(\UNION)
   =& ∑_{u:ℳ} ∏_{z:ℳ} z∈u ↔ ∑_{y : ℳ} y∈x ∧ z ∈ y \\
 τ␣()␣(\UNION)
   =& ∃(u:ℳ)\ ∏_{z:ℳ} z∈u ↔ ∃(y : ℳ)\ y∈x ∧ z ∈ y
\end{align*}

If the structure $(ℳ,∈)$ satisfies the extensionality
axiom, then the property $∏_{z:ℳ} z∈u ↔ ∃(y : ℳ)\ y∈x ∧ z ∈ y$
completely characterises $u$, making
$τ␣()␣(\UNION) ≃ ∑_{u:ℳ} ∏_{z:ℳ} z∈u ↔ ∃(y : ℳ)\ y∈x ∧ z ∈ y$.
However, $∑_{y : ℳ} y∈x ∧ z ∈ y$ is not always implied
by $∃(y : ℳ)\ y∈x ∧ z ∈ y$ so the two axioms remain distinct.

## Axioms of set theory

The axioms of set theory contain a number of axiom schemas, such as
collection or replacement, or (restricted) separation. In set theory
this adds one axiom for each first-order formula. In type theory it is
more convenient to use the higher order features of type theory and
regard these as quantified over all predicates. The result is much
stronger than the original axiom scheme. In the following we will
exploit this extra strength.

#### Definition

Given a structure, $(ℳ,∈)$, and a predicate \newline
$P : ℳ → ℳ → \Type$ and
an element $m : ℳ$,
we define *σ-replacement for $P$ and $m$ in $(ℳ,∈)$* to be the
type \newline
$σ_{P}␣(a)␣((∀x∈a␣∃!y\ P␣x␣y)→∃b␣∀y(y∈b ↔ ∃x∈a\ P␣x␣y))␣m$,
where $σ_{P}$ is $σ$ extended with the clause 
$σ␣Γ␣(P␣x␣y)␣γ ≔ P␣(γ␣x)␣(γ␣y)$, in order to interpret
the $P$ as a predicate symbol.

Define τ-replacement in the same way, substituting τ for σ.

#### Proposition
\label{ac-from-replacement}

If an ∈-structure $(ℳ,∈)$ satisfies extensionality, σ-replacement and has
an ordered pairing operation $〈-,-〉 : ℳ → ℳ → ℳ$, then the following
choice principle, stemming from the so-called *type theoretic choice
principle*^[The type theoretic choice principle is the fact that
 for all $A$, $B$ and $P$, the type
$\left(∏_{a:A}∑_{b:B␣a}P␣a␣b\right) → ∑_{f:∏_{a:A}B␣a}∏_{a:A}P␣a␣(f␣a)$
is inhabited.], holds:

For any $P : ℳ → ℳ → \Type$ and $a,b : ℳ$ if

\begin{align*}
∏_{x : ℳ} x∈a → ∑_{y : ℳ} {y∈B}∧{P␣x␣y}\text{,}
\end{align*}

\noindent then there is a function
$f : ℳ$ with domain $a$ and codomain $b$ such that
$∏_{x : ℳ} P␣x␣(f␣x)$.

*Proof:* Given $a$,$b$ and a proof
$p : ∏_{x : ℳ} x∈a → ∑_{y : ℳ} {y∈B}∧{P␣x␣y}$, apply σ-replacement to $a$
and the predicate $P' : ℳ → ℳ → \Prop$ defined by

\begin{align*}
 P'␣x␣z &≔ ∑_{q : x ∈ a} z =_ℳ {〈x,π₀␣(p␣x␣q)〉}
\end{align*}
The resulting element of $ℳ$ is a function with the desired properties.∎

#### Remark

Proposition \ref{ac-from-replacement} shows us that we cannot hope to
have a model satisfying all axioms of constructive set theory by
interpreting the existential quantifier as Σ-types while at the same
time interpreting equality as the identity type of a mere set. This is
because the above proposition shows that AC will hold, in this
interpretation, and thus if all other CZF axioms hold (in fact
$U$-restricted separation and pairing should suffice), Diaconescu’s
theorem demonstrates that the law of excluded middle holds (for all
$U$-small propositions).

# The model of iterative sets

We recall the definition, from \multisetsref,
of the type of iterated
multisets and the membership relation.

#### Definition {.definition}

\begin{align*}
M ≔ W_UT
\end{align*}
\begin{align*}
&∈_M : M → M → \Set \\
x &∈_M (\sup␣a␣f) ≔ ∑_{i:T␣a} f␣i = x
\end{align*}

The elements of $M$ are the iterative multiset in which another element
of $M$ may occur any $U$-small number of times. This is expressed by
the relation $∈_M$, as follows: Given $x,y : M$ the type $x ∈_M y$ is the
type of occurrences of $x$ in $y$. For instance, if $(x ∈_M y) ≃ 2$ then
$x$ occurs twice in $y$.

What we would like to consider are the elements of $M$ in which every
element occurs at most once. Such a multiset would be called *set-like*.
Being set-like could be expressed in different ways. The most direct
would be to say that $y : M$ is set-like whenever $x ∈_M y$ is a mere
proposition for all other $x : M$. Another way to state this is to say
that the function mapping each instance of an element in $y$ to
the element it represents is an embedding.

The *iterative sets* are those multisets which are hereditarily set-like.
On $M$ we define the predicate $\itset$, recursively, as follows.

#### Definition

\begin{align*}
&\itset : M → \Type \\
&\itset␣(\sup␣a␣f) ≔ \embedding␣f ∧ ∏_{i : T␣a} \itset␣(f␣a)
\end{align*}
For each $x : M$, whenever we have $\itset␣x$, we say that $x$ is an
*iterative set*.

A Σ-type then collects the type of iterative sets as a
subtype of $M$ – with elementhood defined by the restriction of
elementhood for multisets.

#### Definition {.definition}
The type of iterative sets, $V$, is defined as the subtype of $M$ of
iterative sets, or $V ≡ ∑_{x:M} \itset␣x$.

We denote by $∈_V$ the specialisation of $∈_M$ to $V$, that is\newline
$x ∈_V y ≔ π₀␣x ∈_M π₀␣y$.

#### Notation

We will permit ourselves a slight abuse of notation when denoting
elements of $V$. We will for the remainder of this article
use the notation $\sup␣a␣f$ to denote an element of $V$ constructed
by $a : U$ and an embedding $f : T␣a ↪ V$.

#### Remark

Checking that a multiset is an iterative set can be a tedious task,
since it has to be carried out on each level, but the next section will 
give lemmas to make it easier to construct
new iterative sets.

# Basic results

In this section we clarify the basic properties of iterative multisets
and \newline $(V, ∈_V)$. The most important of these is the extensionality of
$V$ with respect to $∈_V$.

We start with reminding ourselves of the extensionality theorem for
multisets.

#### Theorem {.theorem}

$x =_M y ≃ ∏_{z : M} \left( \left(z ∈_M y\right) ≃ \left(z ∈_M y\right)\right)$

*Proof:* See Theorem \ref{multiset-extensionality} in Subsection \ref{multiset-extensionality-section} of \multisetsref.

#### Lemma {.lemma}

The type $\itset␣x$ is a mere proposition for every
$x : M$.

*Proof:* Immediate by induction on $M$ and that $\embedding f$ is a
mere proposition.

#### Lemma {.lemma}

The function $\ap␣π₀ : (x =_V y) → (π₀␣x =_M π₀␣y)$ is an equivalence.

*Proof:* Simple consequence of $\itset$ being a mere proposition.

#### Lemma {.lemma}

For any $x: M$ and $y : V$ we have that the type $x ∈_M π₀y$ is a mere proposition.

*Proof:* Assume $y ≡ (\sup␣a␣f , (p , q))$. Observe that $x ∈_M π₀y$ is
the fibre of $f$ over $y$. Since $p$ proves $f$ to be an embedding, and
embeddings have propositional fibres, then $x ∈_M π₀y$ is a
proposition.∎

#### Lemma {.lemma}
\begin{align*}
\left(∏_{z : V} \left(z ∈_V x ↔ z ∈_V y\right)\right)
≃ \left(∏_{z : M} \left(z ∈_M π₀x ↔ z ∈_M π₀y\right)\right)
\end{align*}

*Proof:* Both sides of the equivalence are mere propositions, so it is
enough to show biimplication. Passing from right to left is trivial, so
we show only the implication from left to right.

Assume $z : M$. If $z ∈ π₀x$ then $z$ must be an iterative set. Thus,
$z ∈ y$ which is to say $z ∈ π₀y$. This demonstrates
$z ∈ π₀x → z ∈ π₀y$. That $z ∈ π₀y → z ∈ π₀x$ is shown symmetrically.
Hence, $z ∈ π₀x ↔ z ∈ π₀y$, which completes the proof.∎

#### Lemma
\label{image-lemma}
Given a small type $a : U$ and a function $f : T␣a → V$,
there is a set ${\image}␣a␣f$ such that

-   for each $i : T␣a$ we have that
    $f␣i ∈ ({\image}␣a␣f)$
-   for any merely propositional predicate $P : V → \Set$, given\hfill\newline
    $∏_{i : T␣a} P␣(f␣i)$ we can prove that
    $∏_{z : V} (z ∈ ({\image}␣a␣f) → P␣z)$.

*Proof:*
 Given $a : U$ and $f : T␣a → V$. Take the image factorisation of $f$, which
can be expressed as a simple higher inductive type. Since
$V$ is locally $U$-small, the image has a code $b : U$. Denote the injection
of the image into $V$ by $g :T b → V$, and define $\image␣a␣f := \sup␣b␣v$.
∎

# V models Myhill’s constructive set theory

In this section we prove that $(V , ∈_V)$ models the axioms 
of Myhill’s Constructive Set Theory (CST), when the existential
quantifiers are interpreted as truncated. In fact, we shall see
that except for a few critical places, positive occurrences of
the existential quantifier can be strengthened to ∑, mostly
because the constructions we make are explicit.

## Extensionality

The lemmas we have proved line up to give the following equivalence:

#### Extensionality

$x =_V y ↔ ∏_{z : V} \left(z ∈ x ↔ z ∈ y\right)$

*Proof:*
\begin{align}
 x =_V y &≃ π₀x = π₀ y \\
         &≃ ∏_{z : M} \left(z ∈_M π₀x ≃ z ∈_M π₀y\right) \\
         &≃ ∏_{z : M} \left(z ∈_M π₀x ↔ z ∈_M π₀y\right) \\
         &≃ ∏_{z : V} \left(z ∈_V π₀x ↔ z ∈_M π₀y\right)
\end{align}∎

## The empty set, natural numbers

The empty set is given by $∅ ≔ \sup 0 e$. The natural
numbers we constructed in \multisetsref, Subsection 
\ref{natural-numbers-multiset} is indeed an iterative set.

## Separation

Restricted separation can be done in $V$ without quotienting, as long as
the separating predicate is merely propositional. Thus, whenever we
separate a formula with existential quantifier in a positive position,
the existence must be unique or truncated. If there are more than one
witness of the statement, the result would be a multiset where the
element occurs once per witness of the statement. The same is true for
disjunction. Notice that if $A + B$ is a proposition then $A$ and $B$
are mutually exclusive (since $l␣a$ and $r␣b$ are always distinct).

#### Proposition

For any $U$-small predicate $P : V i → \Prop$, and $x : V$ 
there is $u : V$ such that for any $z : V$ we have 
$z ∈ u ↔ (P␣z) × (z ∈ x)$.

*Proof:* Assume that $x≡\sup␣A␣f$ and let 
$u ≔ \sup (∑_{a : A} P␣(f␣a)) (f ∘ π₀)$. Since $P∘f$ is a
mere proposition $π₀ : (∑_{a : A} P␣(f␣a)) → A$ is injective,
thus $f ∘ π₀$ is injective.

If $z ∈ u$ then $z = f␣a$ for some $a : A$ such that $P␣(f␣a)$,
and hence $P␣z$ and $z ∈ x$.

If $p:P␣z$ and $q:z ∈ x$ then let $a ≔ π₀ q$,
$((a,p),π₁␣q)$ will prove that $z ∈ u$.∎

### ∈-induction

Induction on $V$ can be performed, even when the predicate is a general
type, not necessarily merely propositional.

#### Proposition

For every predicate $P : V → \Type$, if for each $x : V$
we have that $∏_{y : V} y ∈ x → P␣y$ implies $P␣x$, then we have
$P␣x$ for every $x : V$.

*Proof.* By induction on $V$. Assume $x ≡ \sup␣A␣f$ then by induction
hypothesis $P␣(f␣i)$ for every $i : A$. We must show that if $y ∈ x$ then
$P␣y$. However, $y ∈ x$ means that there is $i$ such that $y = f␣i$ so
we can transport the induction hypothesis to obtain $P␣y$. Thus, we have
shown $∏_{y : V} y ∈ x → P␣y$, which implies $P␣x$ by assumption.∎

## Pairing and Union

In order to show pairing and union we will have to apply the image
construction, previously introduced (Lemma \ref{image-lemma}).
If we applied the constructions of union and pairing for multisets
which we defined in \multisetsref, then the resulting multisets
would not be set-like. Therefore, we need a different construction,
and the most natural is to take quotients to make the multisets back
into sets. The fact that quotienting was not needed for multisets,
may be an indication that iterative multisets is a more natural notion
than iterative sets to consider in type theory.

#### Definition
 Given $x,y : V$, we define
$\{x,y\} ≔ \image␣(1 + 1)␣(\const x + \const y)$.

#### Proposition

For any $x,y : V$ and for each $z : V$ we have that
$z ∈ \{x,y\} ↔ ((x = z) ∨ (y = z))$

*Proof:* Simple consequence of Lemma \ref{image-lemma}.∎

#### Definition
Given $x : V$, where $x ≡ \sup␣A␣f$ define

$∪ x ≔ {\image}␣(∑_{a:A} (\overline{f␣a}))␣
                 (λ〈i,j〉. \widetilde{(f i)} j)$.

#### Proposition

For any $x : V$ and for any $z : V$, we have that
$z ∈ ∪x ↔ ∃y\ ({y ∈ x} ∧ {z ∈ y})$

*Proof:* Simple consequence of Lemma \ref{image-lemma}.∎

#### Remark
The truncated existential quantifier in the above proposition cannot be
strengthened to a Σ-type, since it would mean constructing sections for
almost arbitrary quotients.

## Replacement

#### Proposition
\label{replacement}

For any $a : V$ and $P : V → V → \Prop$, such that
for all $x ∈ a$ there exists a unique $y$ for which $P␣x␣y$,
then there is $b : V$ such that for each $y : V$ we have
that $y ∈ b$ if and only if there exists $x ∈ a$ such that $P␣x␣y$.

*Proof:* Given $a = \sup␣\bar a␣\tilde a$, the assumptions let
us construct a map $T␣\bar a → V$. Using Lemma \ref{image-lemma},
we can construct its image in $V$, which will have the desired
properties.∎

#### Proposition

For any $a : V$ and $F : V → V$, there is $y : V$ such
that for any $z : V$ we have $z ∈ b ↔ ∃(w : V)\ z = F␣w$.

*Proof:* Assume $x ≡ \sup␣a␣f$ and let $b = \image␣(F∘f)$, and apply
Lemma \ref{image-lemma}.∎

## Exponentials

The construction of exponentials in this model is
particularly easy, since a set-theoretical function between two elements,
$(\sup␣A␣f)$ and $(\sup␣B␣g)$, of $V$ boils down to a function in the
type theory $A → B$. Instead of giving a direct
proof, we can lean on \multisetsref, Subsection \ref{multiset-exp}
in order to prove the correctness of this construction.

### Exponentiation

Exponentiation of sets is a special case of exponentiation of multisets.
For multisets, we defined^[
Definition \ref{operation-definition} in
Subsection \ref{multiset-functions} of
\multisetsref]$\operation_{a,b}␣f$ for every $a,b,f : M$,
and showed that there is a multiset $\Exp␣a␣b$, such that\newline
${f ∈ \Exp␣a␣b} ≃\operation_{a,b}␣f$.
We will here show that whenever $a,b,f: V$, we have
$\operation_{π₀␣a,π₀␣b}␣(π₀␣f)$ if an only if $\Fun␣a␣b␣f$, and that in fact
$\itset␣(\Exp␣(π₀␣a)␣(π₀␣b))$, in order to conclude that there is
exponentiation in $V$.

#### Lemma
\label{exp-lemma1}

Whenever $a,b,f: V$, we have that\newline
$\operation_{π₀␣a,π₀␣b}␣(π₀␣f) ↔ \Fun␣a␣b␣f$.

*Proof:* Recall that^[ibid.] for every $a,b,f : M$,

\begin{align*}
\operation_{a␣b}␣f ≔ &\left(∏_z␣z ∈ f → ∑_x␣∑_y␣z=〈 x , y 〉 \right) \\ 
                   ∧ &\left(∏_x␣x ∈ a ≃ ∑_y 〈 x , y 〉 ∈ f \right) \\
                   ∧ &\left(∏_y␣y ∈ b ← ∑_x 〈 x , y 〉 ∈ f \right) 
\end{align*}

If here it said $↔$ instead of $≃$, this would state that $f$ is a total
relation between $a$ and $b$. We thus have to show that this
strengthening is exactly the same as ensuring functionality of a total
relation.

In the middle conjunction, $x∈a$ is a mere proposition since $a$ is an
iterative set by assumption. Thus, $∑_y 〈 x , y 〉 ∈ f$ must
also be a mere proposition, and thus equivalent to
$∃!(y:M) 〈 x , y 〉 ∈ f$. This is exactly the requirement for
a total relation to be functional.∎

*Remark:* Since $\operation_{π₀␣a,π₀␣b}␣(π₀␣f)$ and $\Fun␣a␣b␣f$ are
both mere propositions, the biimplication is also an equivalence of
types.

#### Lemma
\label{exp-lemma2}

For every $a, b : V$ the multiset $\Exp␣(π₀␣a)␣(π₀␣b)$ is an iterative
set.

*Proof:* The definition of $\Exp␣a␣(π₀ b)$ is…

\begin{align*}
\operatorname{Exp}␣(π₀ a)␣(π₀ b) := 
                   \sup␣&(\overline{(π₀␣a)}→ \overline{(π₀␣b)}) \\
                        ␣&(λf. \sup␣\overline{(π₀␣a)} \\
                        ␣&(λi. 〈 \widetilde{(π₀␣a)} i
                              , \widetilde{(π₀␣b)}␣(f␣i)〉))
\end{align*}

Thus, we need to show that the function \newline
$λf. \sup␣\overline{(π₀␣a)}␣(λi. 〈 \widetilde{(π₀ a)} i , \widetilde{(π₀ b)} (f␣i) 〉)$,
which takes a function to its graph, is injective, but this comes down
to that a function is determined by its graph, which follows from
function extensionality.

Next, we need to argue that each graph of each function is itself an
iterative set. Assume that $f : \overline{(π₀␣a)} → \overline{(π₀␣b)}$,
then the function which maps an element $i : \bar (π₀ a)$ to
$〈 \widetilde{(π₀ a)}␣i , \widetilde {(π₀ b)}␣(f i) 〉$ is
injective since $π₀␣a$ is an iterative set. Furthermore,
$〈 \widetilde {(π₀␣a)}␣i , \widetilde{(π₀␣b)}␣(f␣i) 〉$ is an
iterative set since both $a$ and $b$ are iterative sets. Thus, the graph of
$f$ is an iterative set.∎

#### Proposition

For every $a,b : V$ there is a $u : V$ such that for any
$c : V$ there is an biimplication $c ∈ u ↔ \Fun␣a␣b␣c$.

*Proof:* Direct from Lemma \ref{exp-lemma1}, Lemma \ref{exp-lemma2}, and
exponentiation in $M$^[`(M-EXP)` in Subsection
\ref{multiset-exp}, \multisetsref],
letting $u ≔ \Exp␣(π₀ a)␣(π₀ b) ↔ \Fun␣a␣b␣f$.∎

# Collection axioms

In this section we discuss the status of collection axioms of
constructive set theory in our model.

Neither strong collection, nor subset collection seem to hold
in either extreme interpretation of the existential quantifier
(i.e. applying τ or σ).
Interpreting the existential quantifier as $Σ$-types forces us
to make arbitrary choices in an apparently unconstructive way.
Interpreting the existential quantifier as the truncated ∃,
the assumptions become too weak to work with.

We discuss two approaches to solving this problem. On the one
hand, we identify which existential quantifiers need
to be weakened – and which have to remain strict – in
order to get something like the collection axioms to become
provable for our model. On the other hand, we identify axioms
about the type theoretical universe, from which $V$ was constructed,
from which we can derive collection and subset collection in the
truncated form.

## Strong Collection

Together with subset collection, strong collection is the most subtle of
the axioms of constructive set theory. First of all because it is
underspecified: the strong collection axiom states the existence of a
set, but does not define it up to equality.

On an intuitive level, strong collection states that if we have shown
that for all elements of a set $x$ there exists some element with a
certain property, then the proof is in some way an operation which
should have an image in $V$.

The problem is that the proof in first-order logic may not be uniform,
in the sense that the operation may not respect equality of elements.
However, this would not prevent us from taking its image.

In Aczel’s original model one can see this, as his $V$ is a setoid and
we can talk about extensional operations on the underlying type. However,
our $V$ has the identity type as its equality type. This means that if
we interpret the existential quantifier as the Σ-type, then any proof
operation will give rise to an actual function which respects equality.
The image of such an operation is naturally a multiset, and we have seen
how to quotient such a multiset to a set. The only wrinkle of this
approach is that the “strong part” of strong collection has to be
weakened by using a truncated existential quantifier.

### Strong quantifiers

The following proposition is the rendering of strong collection in which
we, as far as our abilities go, interpreted the existential quantifier
as Σ-types.

#### Proposition
\label{collection-balenced}

For any predicate $P : V → V → \Type$ and every $a : V$
if $∏_{x : V} (x ∈ a) → ∑_{y : V} P␣x␣y$ then
there is $b : V$ such that

 #. $∏_{x:V} (x ∈ a) → ∑_{y : V} y ∈ b ∧ (P␣x␣y)$
 #. $∏_{y:V} (y ∈ b) → ∃(x : V)\ x ∈ a ∧ (P␣x␣y)$

*Proof:* Given $a ≡ \sup␣\bar a␣\tilde a$, the assumptions let
us construct a map $T␣\bar a → V$. Using Lemma \ref{image-lemma},
we can construct its image in $V$, which will have the desired
properties.∎

#### Remark

It is well known that, in constructive set theory, collection implies
replacement, and that the converse implication does not hold. But the
strong version of replacement proven for our model in Proposition
\ref{replacement}, formulated for a given ∈-structure, does in fact
entail the collection principle of Proposition
\ref{collection-balenced}, for the same ∈-structure.

### Weak quantifies

@gambino2006  discuss a general approach to
interpreting first-order logic into type theory. For any given
interpretation in their framework, they identify type theoretic
principles corresponding to strong collection and subset collection.
Here we will, in a similar fashion, give a sufficient principle
for our model to satisfy strong collection in the sense of weak
quantifiers.

#### Definition

**Collection Principle for a universe, $U$:**

For every locally $U$-small $B : \Set$,
    given $a : U$ and $P : T␣a → B → \Type$,
   such that $∏_{x:T␣A} ∃(y : B)P␣x␣y$ then
   $∃(r : U) ∃(β : T␣r ↪ B) (∏_{x:T␣a}∃(y:T␣r)P␣x␣(β␣y))∧(∏_{y : T␣r}∃(x : T␣a)P␣x␣y)$.

#### Proposition

The collection principle for $U$ implies strong collection in the following sense:

For any predicate $P : V → V → \Type$ and every $a : V$
if $∏_{x : V} (x ∈ a) → ∃(y : V) P␣x␣y$ then
there merely exists $b : V$ such that

 #. $∏_{x:V} (x ∈ a) → ∃(y : V)\ y ∈ b ∧ (P␣x␣y)$
 #. $∏_{y:V} (y ∈ b) → ∃(x : V)\ x ∈ a ∧ (P␣x␣y)$

*Proof:* Apply the collection principle for $U$ to $P ∘ \bar a$, to
obtain the mere existence $r$ and $β$, which together form $b≔\sup␣r␣β$.
1. and 2. follow from the property of the collection principle for $U$.
∎

## Subset Collection

Subset collection is the principle that for any pair of sets there
is a third set, such that each total relation between the
two has an image in the third. The precise formulation of
the axiom in first-order logic is slightly complicated by
the fact that one cannot quantify over all (or even all definable)
total relations in the middle of a formula. Therefore, the axiom
is an axiom schema quantified over ternary formulas, where the
first parameter is allowed to vary. Our formulation will be
close to the first-order formulation:

\begin{align*}
∀ a,b\ ∃c\ ∀ u\ (∀ x∈ a& → ∃y ∈ b ∧ P␣u␣x␣y)\\
→ (∃d∈c &∧ (∀x∈ a → ∃y ∈ d ∧ P␣u␣x␣y) \\
           &∧ (∀y∈d → ∃x\ x ∈ a ∧ P␣u␣x␣y)
\end{align*}

In the above formula there are three existential quantifiers
to interpret, either as a Σ-type or as a truncated existential.
At one extreme it is possible to give all but the last existential
quantifier the Σ-type interpretation, truncating the last one.

#### Proposition

For every predicate $P : V → V → V → \Type$ and every $a,b : V$
there is $c : V$ such that for all $u : V$ if
$∏_{x : V} (x ∈ a) → ∑_{y : V} P␣u␣x␣y$ then there is $d ∈ c$
such that

 #. $∏_{x:V} (x ∈ a) → ∑_{y : V} y ∈ d ∧ (P␣u␣x␣y)$
 #. $∏_{y:V} (y ∈ d) → ∃(x : V)\ x ∈ a ∧ (P␣u␣x␣y)$

*Proof:* Given $a ≡ \sup␣\bar a␣\tilde a$ and
$b ≡ \sup␣\bar b␣\tilde b$, consider the function
$γ:(\bar a → \bar b) → V$ which maps each $f : \bar a → \bar b$ to
$\image (\tilde b ∘ f)$ in $V$. Let $c ≔ \image␣γ$, which will have the
desired properties.∎

### Weak quantifiers

In the same way we did for strong collection, we identify a principle
of subset collection for our universe which is sufficient to prove
the truncated variation of subset collection for our model.

#### Definition
**Subset Collection Principle for $U$:** 

For each $a,b:U$ there merely exists
   $c:U$ such that for every $P:T␣a→T␣b→\Type$ such that $∏_{x:T␣a}∃(y:T␣b)P␣x␣y$,
   there $∃(r : T␣c) (∏_{x:T␣a}∃(y:T␣r)P␣x␣(α␣y)) ∧ (∏_{y : T␣r}∃(x : T␣a)P␣x␣y)$

#### Proposition

For every predicate $P : V → V → V → \Type$ and every $a,b : V$,
there merely exists a $c : V$ such that for every $u : V$,
if
$∏_{x : V} (x ∈ a) → ∑_{y : V} P␣u␣x␣y$ then there merely exists $d ∈ c$
such that

 #. $∏_{x:V} (x ∈ a) → ∃(y : V)\ y ∈ d ∧ (P␣u␣x␣y)$
 #. $∏_{y:V} (y ∈ d) → ∃(x : V)\ x ∈ a ∧ (P␣u␣x␣y)$

*Proof:* Apply the subset collection principle for $U$ to $a$,$b$ to obtain
the mere existence of $c$, and then use the property of $c$ on the predicate
$P␣u$.

# Equivalence with the HIT-formulation

Section 5.1 of Chapter 10 of the book “Homotopy Type Theory”
\footcite{hottbook}, is dedicated to set theory in the context of
homotopy type theory. The iterative hierarchy is presented there, in the
form of a higher inductive type (HIT). Our $V$ is a HIT-free alternative
to this, and in this section we show that the two types are equivalent.

The HIT formulation bears more than a slight similarity to Aczel’s
original construction of the iterative hierarchy in type theory. Both
are quotients of the type we call $M$. The difference is that while
Aczel uses the untruncated version of the quantifiers and leaves the
quotient a setoid, HTT uses the truncated quantifiers and postulates the
existence of a type with the identity type to match.

In this section we define the bisimulation relation on $M$ of which the
$V$ in HTT Chapter 10 can be seen a quotient. We show that this quotient
is equivalent to our $V$. The proof on this relies on a function,
$\itimage : M → V$, which turns any multiset into an iterative set by
identifying all occurrences of each elements, having first applied the
process inductively on all elements.

#### Definition

We define by induction on $M$ the binary relation:

\begin{align*}
≈ &: M → M → \Type \\
(\sup␣a␣f) ≈ (\sup␣b␣g)
   &≔ \left(∏_{x:T␣a} ∃ (y:T␣b) ({f␣x} ≈ {g␣y})\right) \\
   &× \left(∏_{y:T␣b} ∃ (x:T␣a) ({f␣x} ≈ {g␣y})\right)
\end{align*}
\vspace{3mm}

#### Definition

We define by induction on $M$ the function:

\begin{align*}
 &\itimage : M → V \\
 &\itimage␣(\sup␣a␣f) ≔ \image␣a␣(\itimage ∘ f)
\end{align*}

#### Lemma

For each $x : M$ we have that $x ≈ \itimage␣x$.

*Proof:* For each element of $x$ there is a corresponding element in
$\itimage␣x$, since the type of elements of $\itimage␣x$ is a quotient
of the type of elements in $x$. In the other direction there
just *merely* exists an element of $x$ for each element of $\itimage␣x$,
since it was a quotient. However, this is sufficient to show
that $x ≈ \itimage␣x$.

#### Lemma

For iterative sets $≈$ is equivalent to the identity type. That is, for
every $x,y : M$ given $\itset␣x$ and $\itset␣y$, we have that
$x ≈ y → x =_M y$.

*Proof:* By W-induction. Let $x≡\sup␣a␣f$ and $y≡\sup␣b␣g$, and assume
$x ≈ y$. By extensionality, it suffices to show that for any $z : M$ we
have that $z ∈_M x$ if and only if $z ∈_M y$.

In one direction, if $z ∈_M x$ we know that there is $i : T␣a$ such that
$z =_M f␣i$. Since $z ∈_M y$ is a mere proposition, we can assume from
$x ≈ y$ that there is a $j : T␣b$ such that ${f␣i} ≈ {g␣j}$. By
inductive hypothesis, ${f␣i} ≈ {g␣j} → f␣i =_M g␣j$, and thus
$z =_M g␣j$ which is to say $z ∈_M y$.

The other direction is symmetric in $x$ and $y$.

#### Proposition {#hit-equivalence}

$V ≃ M/≈$

*Proof:* Direct consequence of the previous two lemmas.

#### Remark

That $M/≈$ is equivalent to the HIT formulation in the book
can be seen from Lemma 10.5.5 of “Homotopy Type Theory”\footcite{hottbook}.
Thus, proposition \ref{hit-equivalence} shows that our $V$ is indeed
equivalent to the HIT formulation of $V$.

\printbibliography[heading=subbibliography] 

\end{document}